\documentclass{amsart}
\usepackage{microtype,amssymb,algorithm,algpseudocode,graphicx,color,enumerate,cite}

\def\RR{\ensuremath{\mathbb{R}}}
\newcommand{\CC}{{\mathbb C}}

\newtheorem{theorem}{Theorem}
\numberwithin{theorem}{section} 
\numberwithin{equation}{section}

\newcommand{\R}{{\mathbb{R}}}

\newcommand{\wh}{\widehat}
\newcommand{\f}{\frac}

\newcommand{\lt}{\left}
\newcommand{\rt}{\right}

\newsavebox{\smlmat}
\savebox{\smlmat}{$\left(\begin{smallmatrix}1 & 1 & 1 \\ 0 & 1 & 1 \\ 1 & 3 & 3\end{smallmatrix}\right)$}

\title{Critical points for two-view triangulation}
\author{Hon Leung Lee}
\address{Hon Leung Lee, Mathematics, University of Washington, Seattle, WA 98195}
\email{hllee@uw.edu}

\begin{document}

\maketitle

\begin{abstract}
Two-view triangulation is a problem  of minimizing a quadratic polynomial under an equality constraint.
We derive a polynomial that encodes the local minimizers of this problem using the theory of Lagrange multipliers. 
This offers a simpler derivation of the critical points that are given in Hartley-Sturm \cite{hartley1997triangulation}.
\end{abstract}

\section{Introduction}

Two-view triangulation is the problem of estimating a point $X\in \R^3$ from two noisy image projections; 
see \cite[Chapter 12]{hartley-zisserman-2003} for its significance in structure from motion in computer vision. 
Assuming a Gaussian error distribution, one way to solve the problem is to compute the maximum likelihood estimates (MLE) 
for the true image point correspondences. After that the point $X\in \R^3$ can be recovered via linear algebra 
\cite{hartley-zisserman-2003}. In this paper we study the above problem of finding the MLEs. 
According to the discussion in \cite{aholt-agarwal-thomas} or \cite[Chapter 12]{hartley-zisserman-2003},
the problem is formulated as follows.

Consider a rank two matrix $F\in \R^{3\times 3}$ which is called a {\em fundamental matrix} in multi-view geometry. 
This matrix $F$ encodes a pair of projective cameras \cite[Chapter 9]{hartley-zisserman-2003}. Given two points $u_1,u_2\in \R^2$ which denote the noisy image projections, we solve the problem
\begin{align} \label{triangulation} 
\begin{split}
 \min_{x_1,x_2\in \RR^2} \  & \|x_1 - u_1\|^2_2 + \|x_2- u_2\|^2_2\\
 \text{subject to}   \  \ & \wh{x}_2^\top F \wh{x}_ 1= 0
\end{split}
\end{align}
where $\wh{x}_k := (x_k^\top \  1)^\top\in \R^3$ for $k=1,2$. 
The equation $ \wh{x}_2^\top F \wh{x}_ 1= 0$ is called the {\em epipolar constraint}, which indicates that $x_1$ and $x_2$ are the true image projections under the projective cameras associated with $F$. The minimizers of \eqref{triangulation} are the MLEs for the true image correspondences, assuming the error is Gaussian.

In \cite[Chapter 12]{hartley-zisserman-2003} (or  \cite{hartley1997triangulation})  there is a technique for finding the global minimizers of 
\eqref{triangulation} using a non-iterative approach. They use multi-view geometry to reformulate the problem \eqref{triangulation} as minimizing a fraction in a single real variable say $t$. Using the Fermat rule in elementary calculus, it turns out that the minimizers can be computed via finding the real roots of a polynomial in $t$ of degree 6.

In this note, we view the problem \eqref{triangulation} as minimizing a multivariate quadratic polynomial over one single equality constraint, and then employ the classical method of Lagrange multipliers to locate the potential local minimizers. These candidates are called {\em critical points}. 
For general rank two matrices $F$ and general points $u_1,u_2$, there are six critical points. They can be computed via finding the roots of a polynomial of degree 6 in the Lagrange multiplier.
Assuming that a global minimizer exists, the minimizer of \eqref{triangulation} can be obtained from the critical points.

%

\section{Six critical points for two-view triangulation}

\subsection{Reformulation of the problem \eqref{triangulation}}

Given a fundamental matrix $F\in \R^{3\times 3}$ and $u_1=\begin{pmatrix} u_{11} & u_{12}\end{pmatrix}^\top$, $u_2
= \begin{pmatrix} u_{21} & u_{22}\end{pmatrix}^\top\in \R^2$, 
consider the invertible matrices 
$W_1 :=\lt( \begin{smallmatrix} 1 & 0 & -u_{11} \\ 0 & 1 & -u_{12} \\ 0 & 0 & 1\end{smallmatrix} \rt)$ and 
$W_2 :=\lt( \begin{smallmatrix} 1 & 0 & -u_{21} \\ 0 & 1 & -u_{22} \\ 0 & 0 & 1\end{smallmatrix} \rt)$. 
Note that $\|x_k-u_k\|^2  =  \|\wh{x}_k-\wh{u}_k\|^2$.  
and that problem \eqref{triangulation} is equivalent to the problem
\begin{align*} 
\begin{split}
 \min_{x_1,x_2\in \RR^2} \  & \|W_1 \wh{x}_1\|^2_2 + \|W_2 \wh{x}_2\|^2_2\\
 \text{subject to}   \  \ & \wh{x}_2^\top F\wh{x}_ 1= 0
\end{split}
\end{align*}
For all $k=1,2$, the last coordinate of $W_k\wh{x}_i$ equals one. As a result, we let $y_k\in \R^2$ be such that 
$\wh{y}_k = W_k  \wh{x}_k$. 
Then \eqref{triangulation} 
is further equivalent to the problem
\begin{align} \label{t2} 
\begin{split}
 \min_{y_1,y_2\in \RR^2} \  & \f{1}{2}\lt( \|\wh{y}_1\|^2_2 + \|\wh{y}_2\|^2_2\rt) \\
 \text{subject to}   \  \ & \wh{y}_2^\top F'\wh{y}_ 1= 0
\end{split}
\end{align}
where $F':= W_2^{-\top} F W_1^{-1}=   \lt( \begin{smallmatrix} a & b & c \\ d & e & f\\ g & h & i\end{smallmatrix} \rt)$ is another fundamental matrix.

\subsection{Derivation of a six degree polynomial}

Let $G(y_1,y_2) :=  \f{1}{2}\lt( \|\wh{y}_1\|^2_2 + \|\wh{y}_2\|^2_2\rt)$ and 
$H(y_1,y_2):=\wh{y}_2^\top F'\wh{y}_ 1$. 
The Karush-Kuhn-Tucker (KKT) equation for \eqref{t2} 
is $\nabla G + \lambda \nabla H = 0$ for some $\lambda\in \CC$ called the Lagrange multiplier; see
any nonlinear programming text e.g. \cite{bertsekas1999nonlinear}.
Unwinding this equation we obtain 
a linear system in four variables, namely,
\begin{align} \label{lsls} 
\begin{pmatrix}
1 & 0 & \lambda a & \lambda b \\
0 & 1 & \lambda d & \lambda e \\
\lambda a & \lambda d & 1 & 0 \\
\lambda b & \lambda e & 0 &  1 
\end{pmatrix}
\begin{pmatrix}
y_{21}\\
y_{22}\\
y_{11}\\
y_{12}
\end{pmatrix}
= -\lambda
\begin{pmatrix}
c\\
f\\
g\\
h
\end{pmatrix}
\end{align}
where $y_k = \begin{pmatrix} y_{k1} & y_{k2}\end{pmatrix}^\top$ for $k=1,2$, 
and $\lambda$ is the Lagrange multiplier.
To acquire the critical points we derive a polynomial equation in $\lambda$. It comes from first expressing 
$y_k$, $k=1,2$, in terms of $u_1,u_2,F$ and then substituting these expressions into the epipolar constraint
$ \wh{y}_2^\top F'\wh{y}_ 1= 0$. 
Let $A_\lambda$ be the $4\times 4$ coefficient matrix of the above system. One has  
$$
\det(A_\lambda) = (bd-ae)^2\lambda^4 - (a^2+b^2 + d^2 +e^2) \lambda^2 + 1. 
$$
Define $p_{kl}:= \det(A_\lambda) y_{kl}$ for $k,l=1,2$. By Cramer's rule one has 
\begin{align*}
& p_{21} = \lambda[ (bd-ae)(eg-dh)\lambda^3 + (d^2 c+ e^2 c - adf - bef) \lambda^2 + (ag+bh)\lambda-c]\\
& p_{22} = \lambda[ (bd-ae)(ah-bg)\lambda^3 +  (a^2 f + b^2 f - acd - bce) \lambda^2 + (dg+eh)\lambda-f]\\
& p_{11}= \lambda[ (bd-ae)(ce-bf)\lambda^3 + (b^2 g + e^2 g - abh - deh) \lambda^2 +(ac+df)\lambda-g]\\
& p_{12}= \lambda[ (bd-ae)(af-cd)\lambda^3 + (a^2 h + d^2 h - abg - deg ) \lambda^2 + (bc+ef)\lambda-h ].
\end{align*}
Consider the polynomial 
\begin{align*}
T & := -\det(A_\lambda)^2 \wh{y}_2^\top F' \wh{y}_1 = - p_2^\top F' p_1
\end{align*}
where $p_k := \begin{pmatrix} p_{k1} & p_{k2} & \det(A_\lambda)\end{pmatrix}^\top$ for $k=1,2$. 
Since $\det(A_\lambda)$ is a quartic in $\lambda$, and $p_{kl}$ is also a quartic in $\lambda$ for $k,l = 1,2$, 
we know $T$ is  a polynomial in $\lambda$ of degree at most 8. 
By a careful and slightly tedious computation without using any machines, or by using the following {\tt Macaulay2} \cite{M2} code:

\medskip
\noindent
{\small
{\tt 
R = QQ[a,b,c,d,e,f,g,h,i,L];\\
A = matrix\{\{1,0,L*a,L*b\},\{0,1,L*d,L*e\},\{L*a,L*d,1,0\},\{L*b,L*e,0,1\}\};\\
detA = det A;\\
p21 = det matrix\{\{-L*c,0,L*a,L*b\},\{-L*f,1,L*d,L*e\},\{-L*g,L*d,1,0\},\{-L*h,L*e,0,1\}\};\\
p22 = det matrix\{\{1,-L*c,L*a,L*b\},\{0,-L*f,L*d,L*e\},\{L*a,-L*g,1,0\},\{L*b,-L*h,0,1\}\};\\
p11 = det matrix\{\{1,0,-L*c,L*b\},\{0,1,-L*f,L*e\},\{L*a,L*d,-L*g,0\},\{L*b,L*e,-L*h,1\}\};\\
p12 = det matrix\{\{1,0,L*a,-L*c\},\{0,1,L*d,-L*f\},\{L*a,L*d,1,-L*g\},\{L*b,L*e,0,-L*h\}\};\\
T = -(a*p11*p21+b*p12*p21+c*p21*detA+d*p11*p22+\\
\hspace*{2.5em} e*p12*p22+f*p22*detA+g*p11*detA+h*p12*detA+i*detA*detA);
}
}

\medskip
\noindent
we know the coefficient of 
$\lambda^7$ is zero. The coefficient of $\lambda^8$ is
\begin{align*}
&  -(bd-ae)^2 (eg-dh) (ace-abf+baf - bcd + cbd - cae)  +\\
& -(bd-ae)^2 (ah-bg)(dce - dbf + eaf - ecd + fbd - fae) +\\
 & -(bd-ae)^3 (gce - gbf +haf - hcd + ibd  - iae) =  (bd-ae)^3 \det(F)=0
\end{align*}
since $F$ has rank two.
 This implies $T$ is a polynomial in $\lambda$ of degree at most six. Here we record the explicit expression of $T$: 
\begin{align*}
T = \  &   (bd-ae)^2  (acg + dfg + bch + efh - a^2 i - b^2 i - d^2 i - e^2 i ) \lambda^6 + \\
& a^{2} c^{2} d^{2} \lambda^{5} + c^{2} d^{4}
\lambda^{5} + 2 a b c^{2} d e \lambda^{5} + b^{2} c^{2} e^{2} \lambda^{5} + 2 c^{2} d^{2}
e^{2} \lambda^{5} + c^{2} e^{4} \lambda^{5} - \\
& 2 a^{3} c d f \lambda^{5} - 2 a b^{2} c d f
\lambda^{5} - 2 a c d^{3} f \lambda^{5} - 2 a^{2} b c e f \lambda^{5} - 2 b^{3} c e f
\lambda^{5} - 2 b c d^{2} e f \lambda^{5} - \\
& 2 a c d e^{2} f \lambda^{5} - 2 b c e^{3} f
\lambda^{5} + a^{4} f^{2} \lambda^{5} + 2 a^{2} b^{2} f^{2} \lambda^{5} + b^{4} f^{2}
\lambda^{5} + a^{2} d^{2} f^{2} \lambda^{5} +\\
&  2 a b d e f^{2} \lambda^{5} + b^{2} e^{2}
f^{2} \lambda^{5} + a^{2} b^{2} g^{2} \lambda^{5} + b^{4} g^{2} \lambda^{5} + 2 a b d e
g^{2} \lambda^{5} + 2 b^{2} e^{2} g^{2} \lambda^{5} + \\
& d^{2} e^{2} g^{2} \lambda^{5} +
e^{4} g^{2} \lambda^{5} - 2 a^{3} b g h \lambda^{5} - 2 a b^{3} g h \lambda^{5} - 2 a b
d^{2} g h \lambda^{5} - 2 a^{2} d e g h \lambda^{5} -\\
&  2 b^{2} d e g h \lambda^{5} - 2
d^{3} e g h \lambda^{5} - 2 a b e^{2} g h \lambda^{5} - 2 d e^{3} g h \lambda^{5} + a^{4}
h^{2} \lambda^{5} + a^{2} b^{2} h^{2} \lambda^{5} + \\
& 2 a^{2} d^{2} h^{2} \lambda^{5} +
d^{4} h^{2} \lambda^{5} + 2 a b d e h^{2} \lambda^{5} + d^{2} e^{2} h^{2} \lambda^{5} +
a^{3} c g \lambda^{4} + a b^{2} c g \lambda^{4} +\\ 
& a c d^{2} g \lambda^{4} -5 b c d e g
\lambda^{4} + 6 a c e^{2} g \lambda^{4} + a^{2} d f g \lambda^{4} + 6 b^{2} d f g \lambda^{4} +
d^{3} f g \lambda^{4} - \\
& 5 a b e f g \lambda^{4} + d e^{2} f g \lambda^{4} + a^{2} b c h
\lambda^{4} + b^{3} c h \lambda^{4} + 6 b c d^{2} h \lambda^{4} - 5 a c d e h \lambda^{4} +\\
&  b c
e^{2} h \lambda^{4} - 5 a b d f h \lambda^{4} + 6 a^{2} e f h \lambda^{4} + b^{2} e f h
\lambda^{4} + d^{2} e f h \lambda^{4} + e^{3} f h \lambda^{4} -  a^{4} i \lambda^{4} - \\
& 2 a^{2}
b^{2} i \lambda^{4} -  b^{4} i \lambda^{4} - 2 a^{2} d^{2} i \lambda^{4} - 4 b^{2} d^{2} i
\lambda^{4} -  d^{4} i \lambda^{4} + 4 a b d e i \lambda^{4} - 4 a^{2} e^{2} i \lambda^{4} - \\
& 2
b^{2} e^{2} i \lambda^{4} - 2 d^{2} e^{2} i \lambda^{4} -  e^{4} i \lambda^{4} - 2 c^{2}
d^{2} \lambda^{3} - 2 c^{2} e^{2} \lambda^{3} + 4 a c d f \lambda^{3} + 4 b c e f \lambda^{3} -\\
& 
2 a^{2} f^{2} \lambda^{3} - 2 b^{2} f^{2} \lambda^{3} - 2 b^{2} g^{2} \lambda^{3} - 2
e^{2} g^{2} \lambda^{3} + 4 a b g h \lambda^{3} + 4 d e g h \lambda^{3} - 2 a^{2} h^{2}
\lambda^{3} - \\
& 2 d^{2} h^{2} \lambda^{3} - 3 a c g \lambda^{2} - 3 d f g \lambda^{2} - 3 b c h
\lambda^{2} - 3 e f h \lambda^{2} + 2 a^{2} i \lambda^{2} + 2 b^{2} i \lambda^{2} + \\
& 2 d^{2} i
\lambda^{2} + 2 e^{2} i \lambda^{2} + c^{2} \lambda + f^{2} \lambda + g^{2} \lambda + h^{2} \lambda -  i.
\end{align*}

\subsection{The six critical points}
By solving $T = 0$ for $\lambda$, we get six (complex) solutions (counting multiplicities) for $\lambda$, say $\lambda_1, \ldots, \lambda_6$.
Plugging in these six values of $\lambda$ into the linear system \eqref{lsls}, solving the linear system for $y_1$ and $y_2$, and computing $x_1$ and $x_2$, one obtains the critical points for two-view triangulation. 
If $\det(A_{\lambda_k})\neq 0$ for every $k=1,\ldots,6$ then there are precisely six critical points counting multiplicities.

Now we claim that for general fundamental matrices $F$ and points $u_1,u_2\in \R^2$, there are six distinct critical points for two-view triangulation.
The claim is false if and only if 
the discriminant of $T$ or the resultant 
of $T$ and $\det(A_\lambda)$ are zero polynomials.
Instead of computing the desired discriminant and resultant which depend on $u_1,u_2$ and $F$,  one can find an example of $(u_1,u_2,F)$ such that 
the discriminant of $T$ and the resultant 
of $T$ and $\det(A_\lambda)$ take a nonzero value, that is, 
$\det(A_\lambda)\neq 0$ for every solution $\lambda$ of $T$, and the six critical points obtained are distinct. 
If we consider the data $u_1 = \begin{pmatrix} 0 & 0 \end{pmatrix}^\top$, 
$u_2 = u_1$ and $F= \lt( \begin{smallmatrix} 1 & 1 & 1 \\ 0 & 1 & 1 \\ 1 & 3 & 3\end{smallmatrix}\rt)$, then the polynomial 
$T$ becomes $-2\lambda^6+6\lambda^5+3\lambda^4-12\lambda^3-3\lambda^2+12\lambda-3$, and 
there are six distinct complex critical points for the problem \eqref{triangulation}; see Table \ref{t1}.
\begin{table}
\begin{tabular}{|c|c|c|c|}
\hline
$x_{21}$ & $x_{22}$ & $x_{11}$ & $x_{12}$\\
\hline
$0.0596$ & $-0.0321$ & $-0.312$ & $-0.891$\\
\hline
$-0.0843$ & $-2.06$ & $-0.438$ & $-0.0259$ \\
\hline
$-2.42 + 0.0137i$ & $-1.02 -1.56i$ & $-1.57 + 0.714i$ &  $-1.246 - 1.51i$ \\
\hline
$-2.42 - 0.0137i$ & $-1.02 +1.56i$ & $-1.57 - 0.714i$ &  $-1.246 + 1.51i$ \\
\hline
$-1.69+0.0226i$ & $-0.935 + 0.414i$ & $0.748 + 0.169i$ & $-0.279 - 0.574i$\\
\hline
$-1.69-0.0226i$ & $-0.935 - 0.414i$ & $0.748 - 0.169i$ & $-0.279 + 0.574i$\\
\hline
\end{tabular}
\caption{Six critical points for \eqref{triangulation} when $u_1=(0 \ 0)^\top$, $u_2=u_1$ and
$F=$ \usebox{\smlmat}.} \label{t1}
\end{table}

We summarize the discussion in the following theorem.
\begin{theorem}
For general points $u_1,u_2\in \R^2$ and fundamental matrices $F$, there are six complex critical points
for the problem \eqref{triangulation}.
\end{theorem}

\section{Discussion}

One can make sense of the critical points for $n$-view triangulation where $n$ is greater than two. The authors in 
\cite{stewenius2005hard} (cf. \cite{hartley-kahl}) computed the number of critical points for 2 to 7 view triangulation 
are 6, 47, 148, 336, 638, 1081.  Draisma et al. \cite{DHOST} call this list of numbers the {\em Euclidean distance degrees} of the 
multi-view variety associated to 2 to 7 cameras. They conjecture that the general term of this sequence is
$$
C(n):=\f{9}{2}n^3 - \f{21}{2} n^2 + 8n-4. 
$$
One can apply the B\'{e}zout's theorem to conclude that $C(n)$ has order $n^3$, and our paper verified $C(2)=6$. 
However a proof of the above general formula is still unknown.

\bibliographystyle{plain}
\bibliography{lee}

\begin{thebibliography}{1}

\bibitem{aholt-agarwal-thomas}
C.~Aholt, S.~Agarwal, and R.R. Thomas.
\newblock A {QCQP} approach to triangulation.
\newblock In Andrew~W. Fitzgibbon, Svetlana Lazebnik, Pietro Perona, Yoichi
  Sato, and Cordelia Schmid, editors, {\em ECCV (1)}, volume 7572 of {\em
  Lecture Notes in Computer Science}, pages 654--667. Springer, 2012.

\bibitem{bertsekas1999nonlinear}
D.P. Bertsekas.
\newblock {\em {N}onlinear {P}rogramming}.
\newblock Athena scientific Belmont, 1999.

\bibitem{DHOST}
J.~Draisma, E.~Horobet, G.~Ottaviani, B.~Sturmfels, and R.R. Thomas.
\newblock The {E}uclidean distance degree of an algebraic variety.
\newblock {\em Foundations of Computational Mathematics}, 16:99--149, 2016.

\bibitem{M2}
D.R. Grayson and M.E. Stillman.
\newblock Macaulay2, a software system for research in algebraic geometry.
\newblock Available at http://www.math.uiuc.edu/Macaulay2/.

\bibitem{hartley-zisserman-2003}
R.~Hartley and A.~Zisserman.
\newblock {\em {M}ultiview {G}eometry in {C}omputer {V}ision}.
\newblock Cambridge University Press, second edition, 2003.

\bibitem{hartley1997triangulation}
R.~I. Hartley and P.~Sturm.
\newblock Triangulation.
\newblock {\em Computer vision and image understanding}, 68(2):146--157, 1997.

\bibitem{hartley-kahl}
R.I. Hartley and F.~Kahl.
\newblock Optimal algorithms in multiview geometry.
\newblock In {\em ACCV (1)}, pages 13--34, 2007.

\bibitem{stewenius2005hard}
H.~Stew\'{e}nius, F.~Schaffalitzky, and D.~Nister.
\newblock How hard is 3-view triangulation really?
\newblock In {\em Computer Vision, 2005. ICCV 2005. Tenth IEEE International
  Conference on}, volume~1, pages 686--693. IEEE, 2005.

\end{thebibliography}

\end{document}